\documentclass[12pt,a4paper]{amsart}

\usepackage{amsmath}
\usepackage{amscd}
\numberwithin{equation}{section}

\input{prepictex}
\input{pictex}
\input{postpictex}
\chardef\bslash=`\\ 

\makeatletter
\def\verbatim{\interlinepenalty\@M \@verbatim
  \leftskip\@totalleftmargin\advance\leftskip2pc
  \frenchspacing\@vobeyspaces \@xverbatim}
\makeatother
\newtheorem{theorem}{Theorem}[section]

\newtheorem{proposition}[theorem]{Proposition}

\theoremstyle{definition}
\newtheorem{definition}[theorem]{Definition}

\newtheorem{example}[theorem]{Example}

\newcounter{picture}

\newsymbol\rtimes 226F 

\DeclareMathOperator{\rank}{rank}

\DeclareMathOperator{\coker}{coker}
\DeclareMathOperator{\tor}{tor}

\newcommand{\FF}{{\mathbb F}}

\newcommand{\KK}{{\mathbb K}}
\newcommand{\NN}{{\mathbb N}}

\newcommand{\QQ}{{\mathbb Q}}

\newcommand{\ZZ}{{\mathbb Z}}
\newcommand{\cA}{{\mathcal A}}
\newcommand{\cB}{{\mathcal B}}
\newcommand{\cC}{{\mathcal C}}
\newcommand{\cE}{{\mathcal E}}
\newcommand{\cF}{{\mathcal F}}

\newcommand{\cK}{{\mathcal K}}

\newcommand{\cO}{{\mathcal O}}

\newcommand{\G}{{\Gamma}}
\newcommand{\Om}{{\Omega}}
\newcommand{\g}{{\gamma}}
\newcommand{\s}{{\sigma}}
\newcommand{\w}{{\omega}}

\newcommand{\fp}{{\mathfrak p}}

\newcommand{\abar}{{\overline a}}
\newcommand{\bbar}{{\overline b}}
\newcommand{\cbar}{{\overline c}}

\newcommand{\ow}{{\overline w}}
\newcommand{\oW}{{\overline W}}
\newcommand{\tA}{{\widetilde A}}
\newcommand{\tX}{{\widetilde X}}
\newcommand{\zomatrix}{{$\{0,1\}$-matrix}}

\newcommand{\Ind}{{\bf 1}} 
\newcommand{\id}{{\bf 1}}
\newcommand{\PGL}{{\text{\rm{PGL}}}}

\newcommand{\wt}{\widetilde}
\newcommand{\tv}{\tilde v}

\begin{document}

\title[]{Boundary actions for affine buildings and higher rank Cuntz-Krieger algebras.}

\date{August 29, 1999}
\author{Guyan Robertson }
\address{Mathematics Department, University of Newcastle, Callaghan, NSW 
2308, Australia}
\email{guyan@maths.newcastle.edu.au}
\subjclass{Primary 46L35; secondary 46L55, 22D25, 51E24.}
\keywords{$C^*$-algebra, affine building, Cuntz-Krieger algebra}
\thanks{This research was supported by the Australian Research Council.} 
\thanks{ \hfill Typeset by  \AmS-\LaTeX}

\begin{abstract}
Let $\G$ be a group of type rotating automorphisms of an affine building $\cB$ of type $\wt A_2$. If $\G$
acts freely on the vertices of $\cB$ with finitely many orbits, and if $\Omega$ is the (maximal) boundary of $\cB$, then $C(\Om)\rtimes \G$ is a p.i.s.u.n. $C^*$-algebra. This algebra has a structure theory analogous to that of a simple Cuntz-Krieger algebra and is the motivation for a theory of  higher rank Cuntz-Krieger algebras, which has been developed by T. Steger and G. Robertson. The K-theory of these algebras can be computed explicitly in the rank two case. For the rank two examples of the form $C(\Om)\rtimes \G$ which arise from boundary actions on $\wt A_2$ buildings, the two K-groups coincide.
\end{abstract}      

\maketitle

\section*{Introduction}

Two decades ago J. Cuntz and W. Krieger introduced the class of $C^*$-algebras  which now bears their names \cite{ck}. One reason for the importance of these algebras was their relationship to the classification of topological Markov chains. However they have also proved to be important in several other ways.  Their theory has been refined and extended over the years by many authors from  different points of view. 

The Cuntz-Krieger algebra ${\mathcal O\sb A}$ associated with a nondegenerate $n\times n$ matrix $A$ with entries in $\{0,1\}$
is the universal $C\sp *$-algebra generated by partial isometries $s\sb 1,\cdots,s\sb n$ satisfying 
\begin{subequations}\label{CKrelations}
\begin{eqnarray}
s\sb 1s\sb 1\sp*+\cdots+s\sb ns\sb n\sp * = 1 \label{CKrelationsa}\\
s\sb i\sp *s\sb i = \sum\sb{j=1}\sp nA(i,j)s\sb js\sb j\sp * \label{CKrelationsb}
\end{eqnarray}
\end{subequations}

Cuntz and Krieger proved that ${\mathcal O\sb A}$ is simple if and only if the matrix $A$ is irreducible and not a permutation matrix. It was shown by M. R\o rdam \cite{ror} that simple Cuntz-Krieger algebras are classified up to stable isomorphism by their $K_0$-group. The subsequent classification theorem of E. Kirchberg and C. Phillips \cite{k,k',an} says that purely infinite, simple, separable, unital, nuclear (p.i.s.u.n.) $C^*$-algebras which satisfy the Universal Coefficient Theorem are classified up to isomorphism by their two K-groups together with the class of the identity element in $K_0$. This result applies in particular to simple Cuntz-Krieger algebras.
The K-theory of a Cuntz-Krieger algebra $\mathcal O_A$ can be characterized as follows (see \cite {c3}):
\begin{equation*}
K_0(\mathcal O_A)=(\mbox{finite abelian group})\oplus \ZZ^k ; \:
K_1(\mathcal O_A)=\ZZ^k.
\end{equation*}
The algebras $\mathcal O_A$ are therefore classified up to isomorphism by the group $K_0(\mathcal O_A)$ together with the class of the identity element in $K_0(\mathcal O_A)$.
  
Since p.i.s.u.n. $C^*$-algebras are now relatively well understood, it is of some interest when it happens that such algebras are naturally associated with concrete groups and geometries. Just such a situation has been studied by T. Steger and G. Robertson \cite{rs1,rs2,rs3}. In \cite{rs1} certain group actions on the boundaries of two dimensional buildings were investigated and the corresponding crossed product algebras were seen to be generated by two Cuntz-Krieger subalgebras. Subsequently \cite{rs2} the properties of these geometric examples were abstracted to provide a set of axioms for a class of $C^*$-algebras. These algebras were seen to have a structure theory completely analogous to that of simple Cuntz-Krieger algebras. It is therefore appropriate to refer to them as higher rank Cuntz-Krieger algebras. The development of the theory in \cite{rs2} was closely modeled on the original work of Cuntz and Krieger. The K-theory of rank two Cuntz-Krieger algebras was studied in \cite{rs3}. The results are similar to the those of \cite{c3} and depend upon the fact that a rank two Cuntz-Krieger algebra is stably isomorphic to the crossed product of an AF-algebra by a $\ZZ^2$-action. Suppose that $\G$ is a group of type rotating automorphisms of an affine building $\cB$ of type $\wt A_2$, and that $\G$
acts freely on the vertices of $\cB$ with finitely many orbits.
Suppose that the algebra $\cA$ arises from the boundary action of $\G$.
It follows from symmetry considerations that $K_0(\cA)=K_1(\cA)$. If $\G$ also acts transitively on the vertices of $\cB$ then the class in $K_0(\cA)$ of the identity element has torsion. Extensive computational results are given in \cite{rs3}.

\bigskip

\section{Cuntz--Krieger algebras arising from boundary actions of free groups.}\label{Spielberg}

The geometric construction of higher rank Cuntz-Krieger algebras from group actions on affine buildings was motivated by work of J. Spielberg \cite{s} in the rank one case.
In \cite{s} certain Cuntz-Krieger algebras were exhibited as crossed product algebras arising from actions of  free products of cyclic groups
on totally disconnected spaces. The construction has a particularly simple geometrical interpretation for a free group $\G$ of finite rank acting on its associated tree. The boundary of the tree is a totally disconnected space upon which $\G$ also acts and this action is used to define the relevant crossed product algebra.

\begin{figure}[htbp]
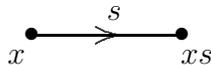
\label{fig0}
\hfil
\centerline{
\beginpicture
\setcoordinatesystem units <1cm, 1cm>  
\setplotarea  x from -1 to 3,  y from -1 to 1
\arrow <10pt> [.2, .67] from  0.8 0  to 1.2 0
\multiput {$\bullet$} at 0 0 *1 2 0 /
\putrule from 0 0 to 2 0
\put {$s$}  [b]   at 1.1  .2
\put {$x$}  [t]   at -.2  -.2
\put {$xs$}  [t]   at 2.2  -.2
\endpicture
}
\hfil
\caption{A labeled edge of the Cayley graph.}
\end{figure}

Consider the specific group $\Gamma = \langle a,b \rangle$, the free group on two generators $a$ and $b$. 
The homogeneous tree $T$ of degree 4 is a Cayley graph of $\Gamma$.  The vertices of $T$ are the elements of $\Gamma$, i.e. reduced words in the generators and their inverses. The edges of $T$ have the form $(x,xs)$, where $x\in \G$ and $s\in 
S = \left \{a, a^{-1}, b, b^{-1} \right \}$.
It is convenient to label the directed edge $(x,xs)$ by the generator $s$ as in Figure \ref{fig0}.

\refstepcounter{picture}

The {\it boundary} $\Omega$ of $T$ can be identified with the set of all infinite reduced
words $\omega = x_1 x_2 x_3 \ldots$, where $x_i \in S$. 
$\Omega$ has a natural compact totally disconnected topology in which a basic open
neighbourhood of $\omega \in \Omega$ consists of those $\omega^\prime \in
\Omega$ whose corresponding infinite word agrees with that of $\omega$ on a finite initial
segment.  Left multiplication by $x \in \Gamma$ defines a homeomorphism of
$\Omega$ and so induces an action $\alpha$ of $\Gamma$ on $C(\Om)$ by
$$
\alpha (x) f(\omega) = f(x^{-1} \omega).
$$

The crossed product $C(\Omega) \rtimes \Gamma$ is the universal $C^*$-algebra  generated by $C(\Omega)$ and the image of a unitary
representation $\pi$ of $\Gamma$, satisfying the relations $\alpha (\g) f = \pi(\g) f \pi(\g)^*$ for
$f \in C(\Omega)$ and $\g \in \Gamma$. It is convenient to write  $\g$ instead of $\pi(\g)$, thereby
identifying elements of $\Gamma$ with unitaries in $C(\Omega) \rtimes \Gamma$. 

\begin{figure}[htbp]
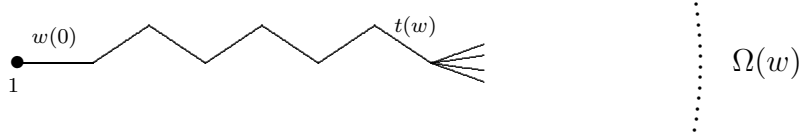
\label{t(w)}
\hfil
\centerline{
\beginpicture
\setcoordinatesystem units <0.5cm, 0.5cm>  
\setplotarea  x from -8 to 10.5,  y from -2 to 3
\put{$_{1}$}   [t] at -8.1 -0.4
\put{$\bullet$}    at -8 0
\put{$_{w(0)}$}   [b] at -7 0.4
\put{$_{t(w)}$}   [b] at 2.6  0.7
\put{$\Omega(w)$}   [l] at 11 0
\putrule from -8 0 to -6 0
\setlinear
\plot   -6 0  -4.5  1   -3  0  -1.5 1   0  0  1.5 1   3  0  4.4 0.5 /
\plot    3  0  4.4 0.2 /
\plot    3  0  4.4 -0.2 /
\plot    3  0  4.4 -0.5 /
\setplotsymbol({$\cdot$}) \plotsymbolspacing=4pt
\circulararc 20 degrees from 10 -1.8 center at 0 0 
\endpicture
}
\hfil
\caption{A basic open subset $\Omega(w)$ of the boundary.}
\end{figure}

\begin{figure}[htbp]
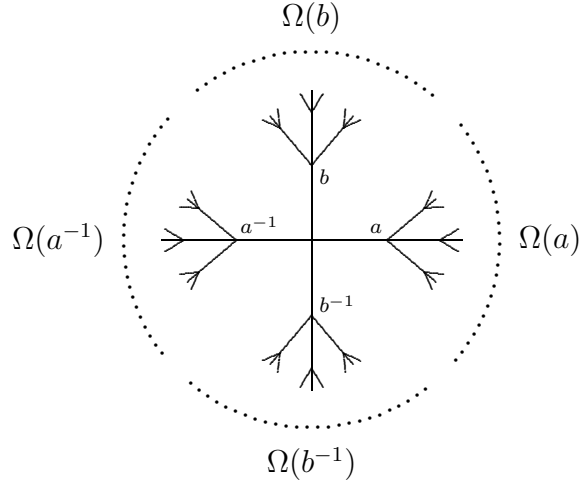
\label{F2graph}
\hfil
\centerline{
\beginpicture
\setcoordinatesystem units <0.5cm, 0.5cm>  
\setplotarea  x from -6 to 6,  y from -6 to 6
\put{${_a}$}   [b,r] at 1.9 0.2
\put{${_b}$}   [t,l] at 0.2 1.9
\put{${_{a^{-1}}}$}   [b,l]  at  -1.9  0.2
\put{${_{b^{-1}}}$}   [b,l]  at   0.2 -1.9
\put{$\Omega(a)$}  [l]  at  5.5 0
\put{$\Omega(b)$}  [b] at  0  5.5
\put{$\Omega(a^{-1})$}  [r] at  -5.5  0
\put{$\Omega(b^{-1})$}  [t] at  0  -5.5
\putrule from -4 0 to 4 0
\putrule from 0 -4 to 0 4
\setlinear
\plot  3.3 1.1  2 0   3.3 -1.1 /
\plot  3.9 0.3  3.4 0   3.9 -0.3 /
\plot   3.2 1.3   3 0.84615  3.5 0.9  /
\plot   3.2 -1.3   3 -0.84615  3.5 -0.9  /
\plot  -3.3 1.1  -2 0   -3.3 -1.1 /
\plot  -3.9 0.3  -3.4 0   -3.9 -0.3 /
\plot   -3.2 1.3   -3 0.84615  -3.5 0.9  /
\plot   -3.2 -1.3   -3 -0.84615  -3.5 -0.9  /
\plot   1.1 3.3  0 2  -1.1   3.3 /
\plot 0.3  3.9  0 3.4  -0.3  3.9  /
\plot  1.3 3.2  0.84615  3  0.9 3.5   /
\plot  -1.3 3.2   -0.84615 3  -0.9 3.5   /
\plot   1.1 -3.3  0 -2  -1.1   -3.3 /
\plot 0.3  -3.9  0 -3.4  -0.3  -3.9  /
\plot  1.3 -3.2  0.84615  -3  0.9 -3.5   /
\plot  -1.3 -3.2   -0.84615 -3  -0.9 -3.5   /
\setplotsymbol({$\cdot$}) \plotsymbolspacing=4pt
\circulararc 80 degrees from 3.83 -3.214 center at 0 0 
\circulararc 80 degrees from 3.214 3.83  center at 0 0
\circulararc 80 degrees from -3.83 3.214 center at 0 0
\circulararc 80 degrees from -3.214 -3.83 center at 0 0
\endpicture
}
\hfil
\caption{The homogeneous tree $T$ of degree four.}
\end{figure}

If $w= w(0)w(1)\dots w(n) \in \G$, where $w(i)\in S$, let $|w|=n$ and let $t(w)=w(n)$, the final letter of the reduced word $w$.
For $w \in \G$ let $\Omega (w)$ be the set of
infinite words beginning with $w$ (Figure \ref{t(w)}).  Then $\Omega(w)$ is open and closed in $\Omega$ and
the sets $\Omega(w)$ for $w \in \G$ form a basis for
the topology of $\Omega$.
The boundary is partitioned into four parts according to the four possible initial letters of $\omega\in \Om$ as shown in Figure \ref{F2graph}.
Denote  by $p_w = \Ind_{\Omega(w)} \in C(\Omega)$ the characteristic function of $\Omega(w)$.

\refstepcounter{picture}

\refstepcounter{picture}

If $u,v \in \G$ and $t(u)=t(v)$, define 
$s_{u,v}=\g p_v \in C(\Omega) \rtimes \Gamma $,
where $\g = uv^{-1}$.
The covariance condition implies that $\g p_v = p_u\g$, so that $s_{u,v}$ is a partial isometry with initial projection $p_v$ and final projection $p_u$.

Let $\cA$ denote the $C^*$-subalgebra of $C(\Omega) \rtimes \Gamma $ generated by 
$\{ s_{u,v}\, ; \, u,v \in \G,\, t(u)=t(v)\}$.
Then $\cA = C(\Omega) \rtimes \Gamma $.  To see this, firstly note that $\cA$ contains $C(\Om)$, since it contains $\{ p_w\ ; \ w \in \G \}$ and this set generates $C(\Om)$ as a $C^*$-algebra.
Also each element $u\in \G$ lies in $\cA$ since 
$$u = \sum_{|x|=|u|+1}up_x = \sum_{|x|=|u|+1}s_{ux,x}.$$

Finally we claim that $\cA$ is a Cuntz-Krieger algebra. For each $x\in S$ let
$$r_x= \sum_{y\in S ;\, |xy|=2}s_{xy,y}= \sum_{y\in S ;\, |xy|=2}xp_y.$$
Then
$$r_xr_x^*=\sum_{y\in S ;\, |xy|=2}p_{xy} = p_x ,$$
$$r_x^*r_x= \sum_{y\in S ;\, |xy|=2}p_y=\sum_{y\in S ;\, |xy|=2}r_yr_y^* .$$
Also
$$\sum_{x\in S}r_xr_x^* = \sum_{x\in S}p_x = \Ind.$$
\noindent Therefore $\{r_x ;\, x\in S\}$ satisfies the classical Cuntz-Krieger relations (\ref{CKrelations}). 

\begin{figure}[htbp]
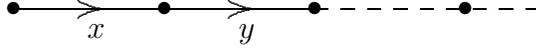
\label{Mone}
\hfil
\centerline{
\beginpicture
\setcoordinatesystem units <1cm, 1cm>  
\setplotarea  x from -4 to 4,  y from -1 to 1
\arrow <10pt> [.2, .67] from  -1.2 0  to  -.8 0
\arrow <10pt> [.2, .67] from  0.8 0  to 1.2 0
\multiput {$\bullet$} at -2 0 *3 2 0 /
\putrule from -2 0 to 2 0
\setdashes
\putrule from 2 0 to 5 0
\put {$x$}  [t]   at -.9 -.2
\put {$y$}  [t]   at 1.1 -.2
\endpicture
}
\hfil
\caption{The condition $M(y,x)=1$.}
\end{figure}

For $u,v\in \G$ with $t(u)=t(v)$, write
 $$r_u=r_{u(0)}r_{u(1)}\dots r_{t(u)}= \sum_{y\in S ;\, |uy|=|u|+1}s_{uy,y}.$$
\noindent Then $r_ur_v^*=s_{u,v}$  . Hence $\cA$ is generated by $\{r_x ; x\in S\}$.  It follows that $C(\Omega) \rtimes \Gamma = \cA = \cO_M$, where 
\begin{equation*}
M=
\begin{pmatrix}
1&0&1&1\\
0&1&1&1\\
1&1&1&0\\
1&1&0&1
\end{pmatrix}.
\end{equation*}
\noindent A geometric interpretation of the condition $M(y,x)=1$  is illustrated by Figure \ref{Mone}.

\refstepcounter{picture}

In the next section higher rank Cuntz-Krieger algebras are defined, based on analogues of the partial isometries $s_{u,v}$ rather than $r_x$.

\bigskip

\section{Higher rank Cuntz-Krieger algebras}\label{HRCK}
 
We begin with some basic notation and terminology from \cite{rs2}.
Let $\ZZ_+$ denote the set of nonnegative integers.
Let $[m,n]$ denote $\{m,m+1, \dots , n\}$, where $m \le n$ are integers. If $m,n \in \ZZ^r$,
say that $m \le n$ if $m_j \le n_j$ for $1 \le j \le r$, and when $m \le n$, let 
$[m,n] = [m_1,n_1] \times \dots \times [m_r,n_r]$. In $\ZZ^r$, let $0$ denote the zero vector
and let $e_j$ denote the $j^{th}$ standard unit basis vector.  We fix a finite set $A$ which we refer to as an ``alphabet''.

A \emph{$\{0,1\}$-matrix} is a matrix with entries in $\{0,1\}$. 
Choose nonzero  $\{0,1\}$-matrices  $M_1, M_2,\dots, M_r$ and denote their elements by
 $M_j(b,a) \in \{0,1\}$ for $a,b \in A$.
If $m,n \in \ZZ^r$ with $m \le n$, let 
$$W_{[m,n]} = \{ w: [m,n] \to A\, ;\ M_j(w(l+e_j),w(l)) = 1\ \text{whenever}\ l,l+e_j\in [m,n] \}.$$
Put $W_m=W_{[0,m]}$ if $m \ge 0$.
Say that an element $w \in W_m$ has \emph{ shape} $m$, and write $\s(w) = m$.
Thus $W_m$ is the set of words of shape $m$, and we identify $A$ with $W_0$ in the
natural way. 
Define the initial and final maps $o: W_m \to A$ and $t: W_m \to A$ by
$o(w) = w(0)$ and $t(w) = w(m)$.

\refstepcounter{picture}
\begin{figure}[htbp]
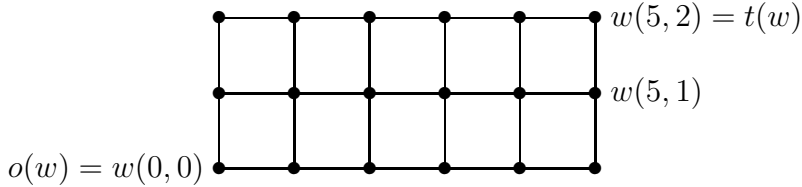
\label{word}
\hfil
\centerline{
\beginpicture
\setcoordinatesystem units <1cm, 1cm>  
\setplotarea  x from -6 to 6,  y from -1 to 1.2
\putrule from -3 -1 to 2 -1
\putrule from -3 0 to 2 0
\putrule from -3 1  to 2  1
\putrule from  -3 -1  to -3 1
\putrule from  -2 -1  to -2 1
\putrule from  -1 -1  to -1 1
\putrule from  0 -1  to -0 1
\putrule from  1 -1  to 1 1
\putrule from  2 -1  to 2 1
\put {$w(5,2)=t(w)$} [l]     at   2.2 1
\put {$w(5,1)$} [l]     at   2.2 0
\put {$o(w)=w(0,0)$} [r]     at   -3.2 -1
\multiput {$\bullet$} at -3 -1 *5  1 0 /
\multiput {$\bullet$} at -3 0 *5  1 0 /
\multiput {$\bullet$} at -3 1 *5  1 0 /
\endpicture
}
\hfil
\caption{Representation of a two dimensional word of shape $m=(5,2)$.}
\end{figure}

Fix a nonempty finite or countable set $D$ (whose elements are ``decorations''),
and a map $\delta : D \to A$. 
Let $\overline W_m = \{ (d,w) \in D \times W_m ;\ o(w) = \delta (d) \}$, the set of
``decorated words'' of shape $m$, and identify
$D$ with $\overline W_0$ via the map $d \mapsto (d,\delta(d))$.
Let $W = \bigcup_m W_m$ and  $\overline W = \bigcup_m \overline W_m$, the sets of all
words and all decorated words respectively.
Define $o: \overline W_m \to D$ and $t: \overline W_m \to A$ by
$o(d,w) = d$ and $t(d,w) = t(w)$. Likewise extend the definition of shape to $\overline W$ by setting
$\s((d,w))=\s(w)$.

Given $j \le k\le l \le m$ and a function $w :[j,m] \to A$, define 
$w \vert _{[k,l]} \in W_{l-k}$ by $w \vert _{[k,l]} = w'$ where
$w'(i) = w(i+k)$ for $0 \le i \le l-k$.
If $\overline w = (d,w) \in \overline W_m$, define
\begin{align*}
\overline w \vert_{[k,l]}& = w\vert_{[k,l]} \in  W_{l-k} \ \text{if}\ k \ne 0, \\ 
\text{and} \qquad 
\ow \vert_{[0,l]}& = (d, w\vert_{[0,l]}) \in  \overline W_l. 
\end{align*}
If $w\in W_l$ and $k\in \ZZ^r$, define $\tau_kw:[k,k+l] \to A$ by 
$(\tau_kw)(k+j)=w(j)$. If $w \in W_l$ where $l \ge 0$ and if $p \ne 0$, say that
$w$ is {\em $p$-periodic} if its $p$-translate,
$\tau_pw$, satisfies $\tau_pw\vert_{[0,l]\cap[p,p+l]} = w\vert_{[0,l]\cap[p,p+l]}$.

Assume that the matrices $M_i$ have been chosen so that the following conditions hold.

\begin{description}
\item[(H0)] Each $M_i$ is a nonzero \zomatrix.

\item[(H1)] Let $u\in W_m$ and $v \in W_n$. If $t(u) =o(v)$ then there exists a unique 
$w\in W_{m+n}$ such that
$$w\vert_{[0,m]}=u \qquad \text{and} \qquad w\vert_{[m,m+n]}=v.$$
\noindent  We write $w=uv$ and say
that the product $uv$ exists. This product is clearly associative. 

\item[(H2)] Consider the directed graph which has a vertex for each $a \in A$
and a directed edge from $a$ to $b$ for each $i$ such that $M_i(b,a) =1$.
This graph is irreducible.

\item[(H3)] Let $p\in \ZZ^r$, $p \ne 0$. There exists some  $w \in W$ which is not $p$-periodic. 
\end{description}
Condition~(H1) holds if the matrices $M_i$, $1\le i\le r$ satisfy the following three conditions \cite[Section~1]{rs2}.
\begin{description}
\item[(H1a)] $M_iM_j=M_jM_i$.

\item[(H1b)] For $i<j$, $M_iM_j$ is a {\zomatrix}.

\item[(H1c)] For $i<j<k$, $M_iM_jM_k$ is a {\zomatrix}.
\end{description}

\begin{definition}\label{CKdefinition}
The $C^*$-algebra $\cA$ is the universal $C^*$-algebra
generated by a family of partial isometries
$\{s_{u,v};\ u,v \in \overline W \ \text{and} \ t(u) = t(v) \}$
satisfying the relations
\begin{subequations}\label{rel1*}
\begin{eqnarray}
{s_{u,v}}^* &=& s_{v,u} \label{rel1a*}\\
s_{u,v}s_{v,w}&=&s_{u,w} \label{rel1b*}\\
s_{u,v}&=&\displaystyle\sum_
{\substack{w\in W;\s(w)=e_j,\\
o(w)=t(u)=t(v)}}
s_{uw,vw} ,\ \text{for} \ 1 \le j \le r   
\label{rel1c*}\\ 
s_{u,u}s_{v,v}&=&0 ,\ \text{for} \ u,v \in \overline W_0, u \ne v. \label{rel1d*}
\end{eqnarray}
\end{subequations}
\end{definition}
  
\noindent The partial isometry  $s_{u,v}$ has initial projection
$s_{v,v}$ and final projection $s_{u,u}$.
If\ $\oW_0=D$ is finite then $\sum_{u\in \oW_0}s_{u,u}$ is an identity for $\cA$ \cite[Section 3]{rs2}.
We refer to $\cA$ as a higher rank Cuntz-Krieger algebra (of rank $r$) despite the fact that if $r=1$ and $D$ is infinite then $\cA$ is in general only stably isomorphic to a classical Cuntz-Krieger algebra.

Suppose that $r=1$, $M=M_1$, $D=A$ and $\delta$ is the identity map. Then the algebra $\cA$ is isomorphic to the simple Cuntz-Krieger algebra $\cO_M$. 
In fact $\cO_M$ is generated by
a set of partial isometries $\{S_a ; a\in A\}$ satisfying the relations 
$S_a^*S_a=\sum_bM(a,b)S_bS_b^*$. If $u\in W$, let
$S_u=S_{u(0)}S_{u(1)}\dots S_{t(u)}$. If $v\in W$ with $t(u)=t(v)$, define  $S_{u,v}=S_uS_v^*$. 
Then the map $s_{u,v}\mapsto S_{u,v}$ establishes an isomorphism of $\cA$ with $\cO_M$.
More generally $\cA$ is isomorphic to a simple Cuntz-Krieger algebra whenever $r=1$ and $D$ is finite. A quick way to see this is to use \cite[Proposition 6.6]{ror}, together with the remarks following Theorem \ref{main1} below.

Tensor products of ordinary Cuntz-Krieger algebras can be identified as higher rank Cuntz-Krieger algebras $\cA$. If $\cA_1$, $\cA_2$ are rank one Cuntz-Krieger algebras, with corresponding irreducible matrices $M_1, M_2$ indexed by alphabets $A_1, A_2$
then $\cA_1 \otimes \cA_2$ is the rank two Cuntz-Krieger algebra $\cA$ arising from
the pair of matrices $M_1\otimes I, I\otimes M_2$ and the alphabet $A_1 \times A_2$.
More interesting examples arise from group actions on affine buildings.
We describe some of these later. 

\begin{theorem}\label{main1} \cite{rs2} The $C^*$-algebra $\cA$ is purely infinite, simple and nuclear. Any nontrivial $C^*$-algebra with
generators $S_{u,v}$ satisfying relations (\ref{rel1*}) is isomorphic to $\cA$.
\end{theorem}

If $D$ is finite, then $\cA$ is unital. 
Therefore $\cA$ is a p.i.s.u.n. $C^*$-algebra and  satisfies the Universal Coefficient Theorem \cite{rs2}.
By the Classification Theorem \cite{k,k'}, $\cA$  is classified by its K-groups and the class of the identity in $K_0$.

Denote by $\cA_D$ the algebra $\cA$ corresponding to a decorating set $D$. Recall that $D$ is finite or countable.
Given any set $D$ of decorations we can obtain another
set of decorations $D \times \NN$, with the decorating map $\delta':D\times \NN
\to A$ defined by $\delta'((d,i)) = \delta(d)$.  It is shown in \cite[Section 5]{rs2} that 
$\cA_{D \times \NN} \cong \cA_D\otimes\cK$. Also, 
for a fixed alphabet $A$ and fixed transition matrices $M_j$, the isomorphism class of 
$\cA_D \otimes \cK$ is independent of $D$.

These facts are used in \cite[Section 6]{rs2} to prove that 
$\cA \otimes\cK \cong \cF \rtimes \ZZ^r$, where $\cF$ is an AF algebra.
The algebra $\cF$ is isomorphic to a subalgebra of $\cA_{A\times \NN}$, and is defined as an inductive limit algebra $\cF = \varinjlim \cC^{(m)}$ where $\cC^{(m)}$ is an isomorphic copy of $\cC = \bigoplus_{a\in A} \cK$, with $\cK$ the compact operators on a separable infinite dimensional
Hilbert space. There is a commuting diagram of inclusions

\begin{equation*}
\begin{CD}
\cC^{(m+e_k)}   @>>>   \cC^{(m+e_j+e_k)}\\
@AAA                                  @AAA\\
\cC^{(m)}                  @>>>      \cC^{(m+e_j)}
\end{CD}
\end{equation*}

\noindent and the action of an element $l \in \ZZ^r$ on $\cF$ maps the subalgebra $\cC^{(m)}$
onto $\cC^{(m+l)}$ for each $m \ge 0$.

\bigskip

\section{A decorated rank one example}\label{decex}

We have seen that one use of decorating sets is to provide a method of passing from $\cA$ to the stabilized algebra $\cA \otimes \cK$. On the other hand nontrivial decorating sets arise even in the rank one case when the construction of Section \ref{Spielberg} is modified to take account of groups of automorphisms of a tree which act freely but not transitively on the vertices of the tree. Here is an example of how such a situation can arise.
 
Let $X$ be a finite connected graph and let $\tX$ be its universal covering graph. Let $\G =\pi (X)$, the fundamental group of $X$. Then $\G$ is a free group which acts freely on $\tX$ with finitely many vertex orbits \cite[Theorem I.9.1]{dd},\cite[Chapter I Section 3]{ser}. Let $\partial \tX$ be the boundary of the tree $\tX$. Then $\G$ acts on $\partial \tX$ and the crossed product algebra $C(\partial \tX) \rtimes \Gamma $ is a rank one Cuntz-Krieger algebra. Let us look at a simple explicit example, where $\G$ is the free group on two generators.

\begin{figure}[htbp]
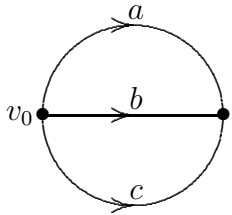
\label{3loop}
\hfil
\centerline{
\beginpicture
\setcoordinatesystem units <0.4cm, 0.4cm>  
\setplotarea  x from -4 to 4,  y from -4 to 4
\put{${a}$}   [b] at 0.1 3.2
\put{${b}$}   [b] at 0.1 0.2
\put{${c}$}   [b] at 0.1 -2.8
\put{$\bullet$}  at -3 0
\put{$v_0$} [r] at -3.3 0
\put{$\bullet$}  at 3 0
\arrow <8pt> [.2, .67] from  -0.1  0 to 0 0
\arrow <8pt> [.2, .67] from  -0.1  3 to 0 3
\arrow <8pt> [.2, .67] from  -0.1  -3 to 0 -3
\putrule from -3 0 to 3 0
\setlinear
\circulararc 360 degrees from 3 0 center at 0 0 
\endpicture
}
\hfil
\caption{The graph $X$.}
\end{figure}

Let $X$ be the directed graph with two vertices and edges $a, b, c$, as illustrated in Figure \ref{3loop}.

\refstepcounter{picture}

\begin{figure}[htbp]
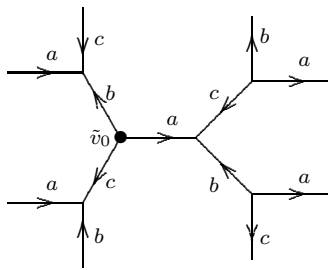
\label{covering}
\hfil
\centerline{
\beginpicture
\setcoordinatesystem units <0.5cm,0.866cm>    
\setplotarea x from -8 to 8, y from -2  to 2         
\putrule from 0 0 to 2 0
\putrule from -3 1 to -1 1
\putrule from -3 -1 to -1 -1
\putrule from 3.5 0.866 to 3.5 1.866
\putrule from 3.5 0.866 to 5.5 0.866
\putrule from 3.5 -0.866 to 3.5 -1.866
\putrule from 3.5 -0.866 to 5.5 -0.866
\putrule from -1 1 to -1 2
\putrule from -1 -1 to -1 -2
\setlinear
\plot -1 1  0 0  -1 -1 /
\plot 3.5 0.866  2 0  3.5 -0.866  /
\put{$_{\tv_0}$}[r] at  -0.2 0
\put{$\bullet$}  at  0 0
\put{$_a$}[b] at  1.4 0.2
\put{$_a$}[b] at  -1.8 1.2
\put{$_a$}[b] at  -1.8 -0.8
\put{$_a$}[b] at  4.9 1.1
\put{$_a$}[b] at  4.9 -0.7
\put{$_b$}[l] at  3.7 1.566
\put{$_c$}[l] at  3.7 -1.566
\put{$_b$}[l] at  -0.7  -1.5
\put{$_c$}[l] at  -0.7  1.5
\put{$_b$}[l] at  -0.4  0.7
\put{$_c$}[l] at  -0.4  -0.7
\put{$_c$}[b] at   2.5  0.6
\put{$_b$}[t] at   2.5  -0.6 
\arrow <6pt> [.3,.67] from  1.2 0 to  1.4 0
\arrow <6pt> [.3,.67] from  -2.0 1 to  -1.8 1
\arrow <6pt> [.3,.67] from  -2.0 -1 to  -1.8 -1
\arrow <6pt> [.3,.67] from  4.7 0.866 to  4.9 0.866
\arrow <6pt> [.3,.67] from  4.7 -0.866 to  4.9 -0.866
\arrow <6pt> [.3,.67] from  3.5 1.466 to 3.5 1.566
\arrow <6pt> [.3,.67] from  3.5 -1.466 to 3.5 -1.566
\arrow <6pt> [.3,.67] from  -1 -1.4 to  -1 -1.3
\arrow <6pt> [.3,.67] from  -1 1.4 to  -1 1.3
\arrow <6pt> [.3,.67] from  -0.6 0.6 to  -0.7 0.7
\arrow <6pt> [.3,.67] from  -0.6 -0.6 to -0.7 -0.7
\arrow <6pt> [.3,.67] from   2.8465 0.5 to  2.6772  0.4
\arrow <6pt> [.3,.67] from  2.8465 -0.5 to  2.6772  -0.4
\endpicture.
}
\hfil
\caption{The universal covering graph $\tX$.}
\end{figure}

Denote by $\abar, \bbar, \cbar$ the opposite edges of $a, b, c$ respectively. The universal covering graph $\tX$ is a homogeneous tree of degree three (Figure \ref{covering}).

\refstepcounter{picture}

By lifting paths in $X$ with initial vertex $v_0$ to paths in $\tX$ with
initial vertex $\tv_0$, the labeling of the edges of $X$ induces a labeling of the edges of $\tX$. The boundary $\partial \tX$ may be identified with the set of all semi-infinite
geodesics in $\tX$ beginning at $\tv_0$. The fundamental group $\G$ of $X$ is the free group on two generators. Choose these generators to be the homotopy classes $[a\bbar], [a\cbar]$.
Any path  in $X$ of even length based at $v_0$  is a loop and its homotopy class in $\G$ is a reduced word in $[a\bbar], [a\cbar]$ and their inverses. For example $[b\cbar]=[a\bbar]^{-1}[a\cbar]$.

The group $\G$ acts by left multiplication on the set of homotopy classes of paths in $X$ beginning at $v_0$. Hence $\G$ acts on $\tX$ and $\partial\tX$.   
The crossed product algebra $C(\partial \tX) \rtimes \Gamma $ is the rank one Cuntz-Krieger algebra $\cA$ constructed in the following way. Let the alphabet be $A=
\{a,\abar,b,\bbar,c,\cbar\}$ and let the decorating set be $D=\{a,b,c\}$ with $\delta : D \to A$ the inclusion map. Let $M$ be the matrix indexed by the elements of $A$, where $M(y,x)=1$ if and only if $xy$ is the labeling of a path of length $2$ in $X$ (that is, a simple loop). Then 
\begin{equation*}
M=
\begin{pmatrix}
0&0&0&1&0&1\\
0&0&1&0&1&0\\
0&1&0&0&0&1\\
1&0&0&0&1&0\\
0&1&0&1&0&0\\
1&0&1&0&0&0
\end{pmatrix}.
\end{equation*}
The set of decorated words $\oW = \{ (d,w) \in D \times W ;\ o(w) = \delta (d) \}$ may clearly be identified with the set of words of the form $w(0)w(1)\dots w(m)$, where $m\in \ZZ_{+}$, $w(0)\in D$, $w(i)\in A,  1\le i\le m$ and $M(w(i+1),w(i))=1$. Thus $\oW$ can be identified with the set of simple paths in $X$ with initial vertex $v_0$. If $w\in \oW$ then $t(w)$ corresponds to the final edge of the path.

The algebra $\cA$ is now defined according to the procedure of Section \ref{HRCK}. 
There is an isomorphism $\phi$ from $\cA$ onto 
$C(\partial \tX) \rtimes \Gamma $ defined  as follows. For $u,v\in \oW$ with $t(u)=t(v)$,
let $[uv^{-1}]$ denote the element of $\G$ represented by the loop in $X$ defined by the path $u$ followed by the inverse path $v^{-1}$. 
Let $\Om(v)\subset \partial \tX$ be the characteristic function of the set of all boundary points which are represented by paths starting at the vertex $\tv_0$ and with initial segment covering  the path $v$. Let $p_v = \Ind_{\Omega(v)} \in C(\partial\tX)$ be the characteristic function of $\Omega(v)$.
Define $\phi(s_{u,v})=[uv^{-1}]p_v$.

In this example and in that of Section \ref{Spielberg} the group $\G$ is the same, but the actions are different. It turns out that the corresponding algebras
have the same K-theory, namely $K_0=K_1=\ZZ^2$, and so the algebras are stably isomorphic. Everything above could also be expressed in the language of groupoids.

\begin{figure}[htbp]
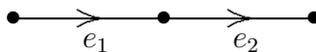
\label{Mmatrix}
\hfil
\centerline{
\beginpicture
\setcoordinatesystem units <1cm, 1cm>  
\setplotarea  x from -4 to 4,  y from -1 to 1
\arrow <10pt> [.2, .67] from  -1.2 0  to  -.8 0
\arrow <10pt> [.2, .67] from  0.8 0  to 1.2 0
\multiput {$\bullet$} at -2 0 *2 2 0 /
\putrule from -2 0 to 2 0
\put {$e_1$}  [t]   at -.9 -.2
\put {$e_2$}  [t]   at 1.1 -.2
\endpicture
}
\hfil
\caption{The geometric condition for $M(\G e_2,\G e_1)=1$.}
\end{figure}

For completeness (and for later comparison) let us state the general one dimensional result.   Assume that $T$ is a tree with fixed base vertex $O$ and boundary $\partial T =\Om$.
Let $\G$ be a group of automorphisms of $T$ that acts freely on the vertex set with finitely many orbits. (In contrast to the example above, $\G$ need not be a free group.) Denote by $\cE$  the set of edges of $T$ and
let  $A= \G\backslash\cE$. Let $D= \{e\in \cE ; O\,  \textrm{ is the initial vertex of}\,  e\}$ and define an injective map $\delta : D \to A$ by $\delta(e)=\G e$.  Define a \zomatrix\, $M$  by $M(b,a)=1$ if and only if $a=\G e_1, b= \G e_2$ for edges $e_1, e_2$ in $T$ lying as indicated in Figure \ref{Mmatrix}. Let $\cA_D$ be the rank one Cuntz-Krieger algebra constructed from these data.

\begin{theorem}\label{treeboundaryaction}
With the above assumptions,  $\cA_D$ is isomorphic to $C(\Om) \rtimes \G $.
\end{theorem}

\refstepcounter{picture}

\bigskip

\section{Affine Buildings}

Let $\cB$ be a finite dimensional simplicial complex, whose maximal simplices we shall call {\em chambers}. All chambers are assumed to have the same dimension and adjacent chambers have a common codimension one face. A {\em gallery} is a sequence of adjacent chambers. $\cB$ is a {\em chamber complex} if any two chambers can be connected by a gallery. $\cB$ is said to be {\em thin} if every codimension one simplex is a face of precisely two chambers. $\cB$ is said to be {\em thick} if every codimension one simplex is a face of at least three chambers. A chamber complex $\cB$ is called a {\em building} if it is the union of a family of subcomplexes, called {\em apartments}, satisfying the following axioms \cite{bro2}.

\begin{description}
\item[(B0)] Each apartment $\Sigma$ is a thin chamber complex with ${\rm dim}\ \Sigma = {\rm dim}\ \cB$.

\item[(B1)] Any two simplices lie in an apartment. 

\item[(B2)] Given apartments $\Sigma$, $\Sigma '$ there exists an isomorphism $\Sigma \to \Sigma '$ fixing $\Sigma \cap \Sigma '$ pointwise.

\item[(B3)] $\cB$ is thick. 
\end{description}

Short readable introductions to the theory of buildings are provided by \cite{bro2,ca,st}.
Detailed introductory texts are \cite{bro1,g}, while \cite{ron} is more advanced.

\begin{proposition}\cite{bro2} The apartments in a building $\cB$ are Coxeter complexes. The Coxeter group is called the Weyl group of $\cB$.
\end{proposition}

\begin{proposition}\cite{bro2} If the apartments are infinite then the building is contractible. The apartments are then affine Coxeter complexes and the building is said to be affine.
\end{proposition}

A building of type $\widetilde A_2$ has apartments which are all Coxeter complexes of type $\widetilde A_2$.
Such a building is therefore a union of two dimensional apartments, each of which may be realized as a tiling of the Euclidean plane by equilateral triangles. 
From now on we shall consider only buildings of type $\widetilde A_2$. These are natural two dimensional analogues of homogeneous trees. In fact a homogeneous tree is  a building of type $\widetilde A_1$. (Such a tree is contractible, its chambers are its edges and the apartments are complete geodesics.) Each vertex $v$ of $\cB$ is labeled with a {\it type} $\tau (v) \in \ZZ/3\ZZ$,
and each chamber has exactly one vertex of each type.
An automorphism $\alpha$ of $\cB$ is said to be {\it type rotating} if there exists $i \in 
\ZZ/3\ZZ$ such that $\tau(\alpha(v)) = \tau(v)+i$ for all vertices $v \in \cB$.

\refstepcounter{picture}
\begin{figure}[htbp]
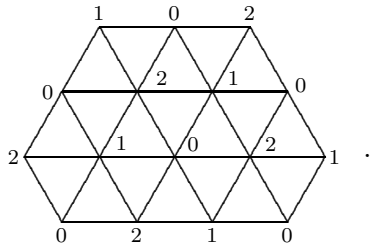
\label{A1}
\hfil
\centerline{
\beginpicture
\setcoordinatesystem units  <0.5cm, 0.866cm>        
\setplotarea x from -2.5 to 5, y from -2 to 2    
\put {$_1$} [b] at -2  2.1
\put {$_0$} [b] at  0  2.1
\put {$_2$} [b] at  2  2.1
\put {$_0$} [l] at 3.2 1.1
\put {$_1$} [l] at 4.1 0
\put {$_0$} [t] at 3  -1.1
\put {$_1$}  [t] at 1  -1.1
\put {$_2$} [t] at  -1  -1.1
\put {$_0$} [t] at -3 -1.1 
\put {$_2$} [r]  at  -4.1 0
\put {$_0$} [r] at -3.2  1
\put {$_2$} [l,b] at -0.5 1.1
\put {$_1$}  [l,b] at 1.4  1.1
\put {$_2$} [l] at  2.4 0.2
\put {$_1$}  [b]  at  -1.4 0.1
\put {$_0$} [t]  at  0.5  0.3
\putrule from -2   2     to  2  2
\putrule from  -3 1  to 3 1
\putrule from -4 0  to 4 0
\putrule from -3 -1  to  3 -1
\setlinear
\plot  -4 0  -2 2 /
\plot -3 -1  0 2 /
\plot -1 -1  2 2 /
\plot  1 -1  3 1 /
\plot  3 -1  4 0 /
\plot  -4 0  -3 -1 /
\plot  -3 1   -1 -1 /
\plot  -2 2   1 -1 /
\plot  0 2   3 -1 /
\plot  2 2   4 0 /
\endpicture.
}
\hfil
\caption{Part of an apartment showing vertex types.}
\end{figure}

A \emph{sector} (or \emph{Weyl chamber}) is a
$\frac{\pi}{3}$-angled sector made up of chambers in some apartment (Figure \ref{A2}).
Two sectors are  \emph{equivalent} (or parallel) if their
intersection contains a sector. (In a tree, sectors are semi-infinite geodesics.)
\refstepcounter{picture}
\begin{figure}[htbp]\label{A2}
\hfil
\centerline{
\beginpicture
\setcoordinatesystem units <0.5cm,0.866cm>   
\setplotarea x from -4 to 4, y from  0 to  3.5  
\putrule from -1 1 to 1 1
\putrule from -2 2 to 2 2
\setlinear \plot -2 2  0 0  2 2  /
\setlinear \plot -1 1  0 2  1 1 /
\setdashes
\setlinear \plot -3.5 3.5  -2 2  -1 3  0 2  1 3  2 2  3.5 3.5 /
\putrule from  -3 3  to  3 3
\endpicture.
}
\hfil
\caption{A sector in a building $\cB$ of type $\tA_2$.}
\end{figure}

The boundary $\Omega$ of $\cB$ is defined to be the set of equivalence classes of sectors in $\cB$. 
In $\cB$ fix some vertex $O$.
For any $\omega \in \Omega$ there is a unique sector $[O,\omega)$ in the
class $\omega$ having base vertex $O$ \cite[Theorem 9.6]{ron}.
The boundary $\Omega$ is a totally disconnected compact Hausdorff space with a base for the
topology given by sets of the form
$$
\Omega(v) = \left \{ \omega \in \Omega : [O,\omega) \ \text{ contains } v \right \}
$$
where $v$ is a vertex of $\cB$ \cite[Section 2]{cms}. 
 
\refstepcounter{picture}
\begin{figure}[htbp]
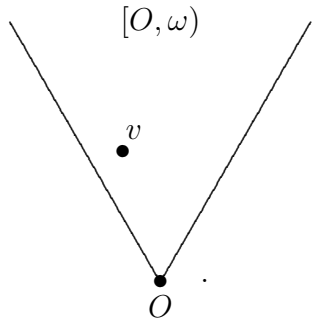
\label{A3}
\hfil
\centerline{
\beginpicture
\setcoordinatesystem units <0.5cm, 0.866cm>
\put {$\bullet$} at 0 0
\put {$\bullet$} at -1 2
\put {$[O,\omega)$}  at 0 4
\put {$O$} [t] at 0 -0.2
\put {$v$}  at  -0.7 2.3
\setlinear \plot -4 4  0 0  4 4 /
\endpicture.
}
\hfil
\caption{The sector $[O,\omega)$, where $\omega\in \Om(v)$.}
\end{figure}

\bigskip

\section{$\tA_2$ groups}\label{A2tildegroups}

Suppose that $\cB$ is a building of type $\tA_2$ and that  $\G$ is a group of type rotating automorphisms of $\cB$ which acts freely and transitively on the vertex set of $\cB$. Such groups $\G$ are called  $\tA_2$ groups. They are good candidates to be rank two analogues of finitely generated free groups, which act in a similar way on buildings of type $\tA_1$ (trees): an automorphism of a tree is automatically type rotating.  The theory of $\tA_2$ groups has been developed in detail in \cite{cmsz,cms}. The $\tA_2$ groups have a detailed combinatorial structure which makes them an ideal place to attack problems involving higher rank groups. For example \cite{cms} proved that $\tA_2$ groups have Kazhdan's Property (T) and obtained exact Kazhdan constants, without the use of an embedding in any linear group. The $\tA_2$ groups were the first examples of higher rank groups known to have property (RD) \cite{rrs}. Also $\tA_2$ groups are a natural place to prove higher rank analogues of results for von Neumann algebras associated with free groups \cite{rr,rs}.

\begin{figure}[htbp]
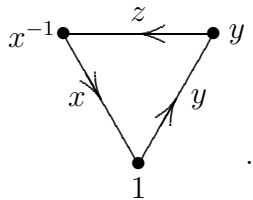
\label{A4}
\hfil
\centerline{
\beginpicture
\setcoordinatesystem units <1cm, 1.732cm>
\put {$\bullet$} at 0 0
\multiput {$\bullet$} at -1 1 *1 2 0 /
\put {$1$} [t] at 0 -0.1
\put {${x^{-1}}$} [r] at -1.1 1 
\put {$y$} [l] at 1.2 1
\arrow <10pt> [.2, .67] from  0.2 1 to 0 1
\arrow <10pt> [.2, .67] from  -0.7 0.7 to -0.5 0.5
\arrow <10pt> [.2, .67] from  0.3 0.3 to 0.5 0.5
\put {$x$} [r ] at -0.7 0.5
\put {$y$} [ l] at 0.7 0.5
\put {$z$} [ b] at 0 1.1
\putrule from -1 1 to 1 1
\setlinear \plot -1 1 0 0 1 1 /
\endpicture.
}
\hfil
\caption{A chamber based at $1$.}
\end{figure}

The $1$-skeleton of  the building $\cB$ is the Cayley graph of the $\tA_2$  group $\G$ with respect to a canonical set $P$ of  $(q^2+q+1)$ generators, where $q$ is a prime power. The set $P$ may be identified  with the set of points of a finite projective plane $(P,L)$ of order $q$.
There are $q^2 + q + 1$ points (elements of $P$) and $q^2+q+1$ lines (elements of $L$).  Each point lies on $q+1$ lines and each line contains $q+1$ points. 
The set of lines $L$ is identified with $\{x^{-1};\ x\in P\}$. The relations satisfied by the elements of $P$ are of the form $xyz=1$. There is such a relation if and only if $y\in x^{-1}$:\; that is the point $y$ is incident with the line $x^{-1}$ in the projective plane $(P,L)$. See Figure \ref{A4}, which illustrates a typical chamber based at the identity element $1\in\cB$. As usual vertices are identified with elements of $\G$ and a directed edge of the form $(a,as)$ with $a\in\G$ is labeled by a generator $s\in P$.  
\refstepcounter{picture}

If $q=2$ there are eight $\tA_2$ groups $\G$, all of which embed as lattices in a linear group $\PGL (3,\FF)$ over a local field $\FF$. If $q=3$ there are 89 possible $\tA_2$ groups, of which 65  have buildings which are not associated with linear groups \cite{cmsz}.

\begin{figure}[htbp]
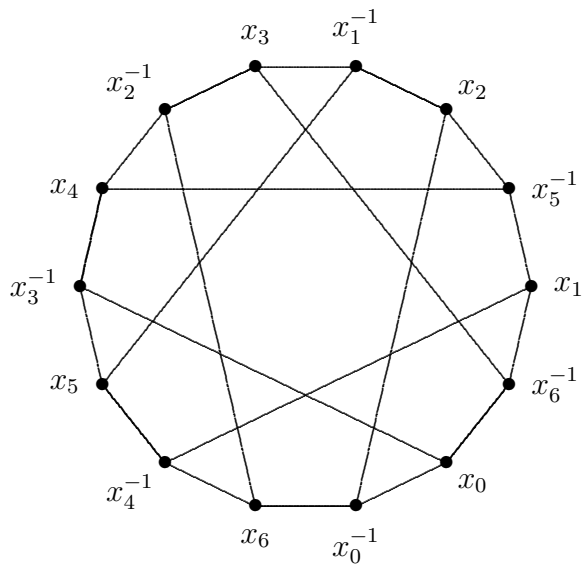
\label{A5}
\hfil
\centerline{
\beginpicture
\setcoordinatesystem units <3cm,3cm> 
\setplotarea  x  from  -1.5  to  1.5, y  from  -1.1  to  1.1
\put {$\bullet$} at 1 0
\put {$\bullet$} at 0.9009690323808488999   0.4338833975744138177 
\put {$\bullet$} at 0.6234903946185663099   0.781831009757469557 
\put {$\bullet$} at 0.222522042695637753     0.9749276591186449509 
\put {$\bullet$} at -0.2225194556367689247  0.9749282495974337219 
\put {$\bullet$} at -0.6234883199575836131   0.7818326642424642541
\put {$\bullet$} at -0.9009678810291215377   0.433885788375114562 
\put {$\bullet$} at -1 0
\put {$\bullet$} at 0.9009690323808488999   -0.4338833975744138177 
\put {$\bullet$} at 0.6234903946185663099   -0.7818310097574695578 
\put {$\bullet$} at 0.222522042695637753     -0.9749276591186449509 
\put {$\bullet$} at -0.2225194556367689247  -0.9749282495974337219 
\put {$\bullet$} at -0.6234883199575836131   -0.781832664242464254
\put {$\bullet$} at -0.9009678810291215377   -0.433885788375114562 
\put{$x_1$} [l] at 1.1 0
\put{$x_5^{-1}$} [l] at 1.0009690323808488999   0.4338833975744138177 
\put{$x_2$} [l,b] at 0.6734903946185663099   0.8318310097574695578 
\put{$x_1^{-1}$} [[b] at 0.222522042695637753     1.0749276591186449509 
\put{$x_3$} [b] at -0.2225194556367689247  1.0749282495974337219 
\put{$x_2^{-1}$} [b,r] at -0.6734883199575836131   0.8318326642424642541
\put{$x_4$} [r] at -1.0009678810291215377   0.433885788375114562 
\put{$x_3^{-1}$} [r] at -1.1 0
\put{$x_6^{-1}$} [l] at 1.0009690323808488999   -0.4338833975744138177 
\put{$x_0$} [l,t] at 0.6734903946185663099   -0.8318310097574695578 
\put{$x_0^{-1}$} [t] at 0.222522042695637753     -1.0749276591186449509 
\put{$x_6$} [t] at -0.2225194556367689247  -1.0749282495974337219 
\put{$x_4^{-1}$} [t,r] at -0.6734883199575836131   -0.8318326642424642541
\put{$x_5$} [r] at -1.0009678810291215377   -0.433885788375114562 
\setlinear \plot  1 0 
 0.9009690323808488999   0.433883397574413817    
 -0.9009678810291215377   0.433885788375114562   
 -1 0    
0.6234903946185663099   -0.7818310097574695578   
 0.9009690323808488999   -0.4338833975744138177  
 -0.2225194556367689247  0.9749282495974337219   
 -0.6234883199575836131   0.7818326642424642541  
-0.2225194556367689247  -0.9749282495974337219   
 0.222522042695637753     -0.9749276591186449509 
 0.6234903946185663099   0.7818310097574695578   
 0.222522042695637753     0.9749276591186449509  
 -0.9009678810291215377   -0.433885788375114562   
-0.6234883199575836131   -0.7818326642424642541  
1 0 
 0.9009690323808488999   -0.4338833975744138177  
0.6234903946185663099   -0.7818310097574695578   
0.222522042695637753     -0.9749276591186449509 
-0.2225194556367689247  -0.9749282495974337219   
-0.6234883199575836131   -0.7818326642424642541  
-0.9009678810291215377   -0.433885788375114562   
 -1 0    
 -0.9009678810291215377   0.433885788375114562   
-0.6234883199575836131   0.7818326642424642541  
-0.2225194556367689247  0.9749282495974337219   
 0.222522042695637753     0.9749276591186449509  
0.6234903946185663099   0.7818310097574695578   
 0.9009690323808488999   0.433883397574413817   /  
\endpicture
}
\hfil
\caption{The projective plane of nearest neighbours of $1$ for the group $\G ({\rm C.1})$.}
\end{figure}

\begin{example}\label{C1} The group C.1 of \cite{cmsz}, which we shall denote $\G ({\rm C.1})$ has presentation
$$\langle
x_i, 0\le i \le 6\,
|\,
x_0x_0x_6,
x_0x_2x_3,
x_1x_2x_6,
x_1x_3x_5,
x_1x_5x_4,
x_2x_4x_5,
x_3x_4x_6
\rangle.$$
For this group $q=2$, and there are $q^2+q+1=7$ generators. We have chosen this group as an example having the smallest possible number of generators (7) and relations (7).  Moreover $\G ({\rm C.1})$ is a lattice subgroup of $G=\PGL (3,\QQ_2)$, where $\QQ_2$ is the field of $2$-adic numbers \cite{cmsz}. 
The vertices of $\cB$ can be identified with the discrete space $G/K$, where 
$K=\PGL (3,\ZZ_2)$, and $\ZZ_2$ is the ring of $2$-adic integers. The boundary $\Om$ of $\cB$ can be identified with $G/B$ where $B$ is a minimal parabolic subgroup of $G$. Thus $G$ acts naturally on the building $\cB$ and its boundary.
A detailed exposition of these facts from an analyst's point of view is given in \cite{st}.   

Figure \ref{A5} illustrates the set of all fourteen neighbours of $1$ in the Cayley graph of $\G ({\rm C.1})$. The fact that there is an edge between $x_2$ and $x_0^{-1}$, for example, is a consequence of the relation $x_0x_2x_3=1$, that is $x_2x_3=x_0^{-1}$.  There are 21 edges in total, each lying in precisely one of the 21 chambers in $\cB$ which contain the vertex $1$. It is worth noting that although we have focused on the vertex $1$, the set of nearest neighbours of any vertex in $\cB$ also has the same structure of a finite projective plane.
\end{example}

\refstepcounter{picture}

\bigskip

\section{Algebras arising from boundary actions on $\tA_2$ buildings}\label{A2boundary}

We are now in a position to describe the class of rank two Cuntz-Krieger algebras which provided the motivation for the general theory of \cite{rs2}.
Theorem \ref{mainboundaryaction} is a rank two version of the examples of Sections \ref{Spielberg} and \ref{decex}.

\begin{theorem}\label{mainboundaryaction} \cite[Theorem 7.7]{rs2}
Let $\cB$ be a building of type $\tA_2$ with boundary $\Om$.
Let $\G$ be a group of type rotating automorphisms of $\cB$ that acts freely on the vertex set with finitely many orbits. Then there is an alphabet $A$, a decorating set $D$ and matrices $M_1, M_2$ such that conditions {\rm(H0-H3)} are satisfied and the corresponding rank two Cuntz-Krieger algebra $\cA$ is isomorphic to $C(\Om) \rtimes \G $.
\end{theorem}

For simplicity, consider the case where the action of $\G$ is also transitive on the vertex set, that is where $\G$ is an $\tA_2$ group, and the $1$-skeleton of $\cB$ is the Cayley graph of $\G$ relative to the generating set $P$. This is a two dimensional analogue of the situation described in Section \ref{Spielberg} and the decorating set is trivial, i.e. $D=A$ . For full generality, with a free but not necessarily transitive action of $\G$ on the vertices and a nontrivial decorating set, see \cite{rs2}; compare also with Section \ref{decex}.

\begin{figure}[htbp]
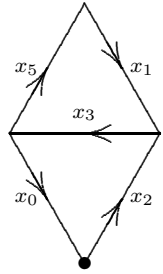
\label{A6}
\hfil
\centerline{
\beginpicture
\setcoordinatesystem units <1cm, 1.732cm>
\setplotarea  x from -2 to 2,  y from 0 to 2
\put {$_{x_3}$} [b ] at 0 1.1
\put {$_{x_2}$} [l ] at 0.6  0.5
\put {$_{x_0}$} [ r] at -0.6  0.5
\put {$_{x_1}$} [l ] at 0.6  1.5
\put {$_{x_5}$} [ r] at -0.6  1.5
\put {$\bullet$}  at 0 0
\arrow <10pt> [.2, .67] from  0.2 1 to 0 1
\arrow <10pt> [.2, .67] from  -0.7 0.7 to -0.5 0.5
\arrow <10pt> [.2, .67] from  0.3 0.3 to 0.5 0.5
\arrow <10pt> [.2, .67] from  0.3 1.7 to 0.5 1.5
\arrow <10pt> [.2, .67] from  -0.7 1.3 to -0.5 1.5
\putrule from -1 1 to 1 1 
\setlinear \plot -1 1 0 0 1 1  0 2 -1 1 /
\endpicture
}
\hfil
\caption{A tile $a\in A$ for the group $\G ({\rm C.1})$.}
\end{figure}

Identify elements of $\G$ with vertices of the building $\cB$. The alphabet $A$ is defined to be the set of $\G$-equivalence classes of basepointed parallelograms in $\cB$,  as illustrated in Figure \ref{A6}. We refer to such an element of $A$ as a tile.
Each tile has a unique representative labelled parallelogram based at a fixed vertex as in Figure \ref{A6}, where each edge label is a generator of $\G$.
 The combinatorics of the finite projective plane $(P,L)$ shows that there are precisely $q(q+1)(q^2+q+1)$ tiles $a\in A$.

\refstepcounter{picture}

\begin{figure}[htbp]
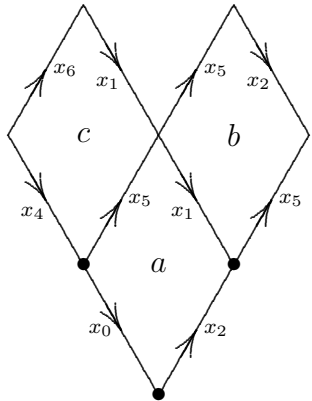
\label{A7}
\hfil
\centerline{
\beginpicture
\setcoordinatesystem units <1cm, 1.732cm>
\setplotarea  x from -2 to 2,  y from -1 to 2
\put {$_{x_0}$} [r ] at -0.6 -0.5
\put {$_{x_2}$} [ l] at 0.6 -0.5
\put {$_{x_1}$} [r ] at 0.5  0.4
\put {$_{x_5}$} [ l] at 1.6  0.5
\put {$_{x_4}$} [r ] at -1.5  0.4
\put {$_{x_5}$} [ l] at -0.4  0.5
\put {$_{x_6}$} [ l] at -1.4  1.5
\put {$_{x_1}$} [r ] at -0.5  1.4
\put {$_{x_5}$} [ l] at 0.6  1.5
\put {$_{x_2}$} [r ] at 1.5  1.4
\put {$a$}  at   0  0
\put {$b$}  at   1  1
\put {$c$}  at   -1  1
\put {$\bullet$}  at 0 -1
\put {$\bullet$}  at -1 0
\put {$\bullet$}  at 1 0
\arrow <10pt> [.2, .67] from  -0.7 -.3 to -0.5 -0.5
\arrow <10pt> [.2, .67] from  0.3 -.7 to 0.5 -0.5
\arrow <10pt> [.2, .67] from  0.3 .7 to 0.5 .5
\arrow <10pt> [.2, .67] from  -0.7 .3 to -0.5 .5 
\arrow <10pt> [.2, .67] from  1.3 1.7 to 1.5 1.5
\arrow <10pt> [.2, .67] from  -1.7 1.3 to -1.5 1.5
\arrow <10pt> [.2, .67] from  -1.7 0.7 to -1.5 0.5 
\arrow <10pt> [.2, .67] from  -0.7 1.7 to -0.5 1.5 
\arrow <10pt> [.2, .67] from  1.3 0.3 to 1.5 0.5
\arrow <10pt> [.2, .67] from  0.3 1.3 to 0.5 1.5
\setlinear \plot -1 0  0 -1   1 0  0 1  -1 0 /
\setlinear \plot -1 0  -2 1  -1 2  0 1 /
\setlinear \plot 1 0  2 1  1 2  0 1 /
\endpicture
}
\hfil
\caption{ $M_1(b,a)=1$, $M_2(c,a)=1$.}
\end{figure}

Suppose that $\G$ is the group $\G ({\rm C.1})$ of Example \ref{C1}. Then $q=2$ and $|A|=42$. 
The transition matrices $M_1$, $M_2$ are defined as follows. If $a,b\in A$ we have $M_1(b,a)=1$ if and only if there are labeled parallelograms representing $a,b$ in the building $\cB$ which lie as shown in Figure \ref{A7}. In that diagram we have chosen edge labels representing specific choices for $a, b$. If no such diagram is possible then $M_1(b,a)=0$. Figure \ref{A7} also illustrates the case $M_2(c,a)=1$. Examination of the edge labels shows that we also have $M_2(b,a)=1$ but that $M_1(c,a)=0$. This geometric definition of the transition matrices is the exact analogue of the one dimensional situation described by Figure \ref{Mone}. 

\refstepcounter{picture}

\begin{figure}[htbp]
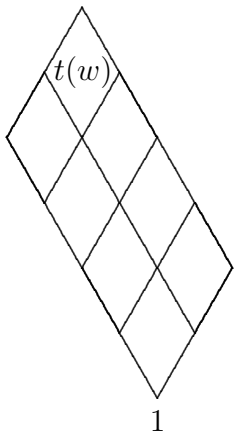
\label{parallelogram}
\hfil
\centerline{
\beginpicture
\setcoordinatesystem units <0.5cm,0.866cm>    
\put{$1$}  [t] at  1 -2.2
\setplotarea x from -6 to 6, y from -2 to 4         
\setlinear
\put{$t(w)$}  [r] at  -.2  3
\plot 1  -2    3  0   2 1   0 -1   -1 0  1 2  0 3  -2 1  -3 2  -1 4  3 0
  2 -1  -2 3  -3  2  1  -2  /
\endpicture
}
\hfil
\caption{A parallelogram $\fp$ and the terminal letter of the word $w=w(\fp)$.}
\end{figure}

Let $\fp$ be a parallelogram based at $1$ in some apartment of $\cB$. Then $\fp$ is a union of parallelograms representing tiles from the alphabet $A$ (Figure \ref{parallelogram}). Associated to $\fp$ there is a two dimensional word $w=w(\fp)$,  as in Section \ref{HRCK}. The map $\fp \mapsto w(\fp)$ is bijective, and by abuse of notation we identify $\fp$ with $w(\fp)$. For example in Figure \ref{A7}, the two letters $a,b$ define a word $w\in W_{(1,0)}$, with $w(0,0)=a$ and $w(1,0)=b$, whereas the two letters $a,c$ define a word $w\in W_{(0,1)}$. If $w=w(\fp)\in W$, let the terminal letter $t(w)\in A$ be the tile of the  parallelogram $\fp$ farthest from $1$ (Figure \ref{parallelogram}).
Also let $\Om(w)=\{ \w\in \Om\, ; \, \fp\subset [1,\w)\}$, the set of boundary points whose representative sectors based at $1$ contain $\fp$.

\refstepcounter{picture}

We can now describe the isomorphism $\phi: \cA \to C(\Om)\rtimes \G$.  If $w_1, w_2 \in W$ with $t(w_1)=t(w_2)=a\in A$, let $\g \in \G$ be the unique element such that $\g t(w_1)=t(w_2)$. Then
\begin{equation}\label{isomorphism} 
\phi(s_{w_2,w_1})=\g \Ind _{\Om(w_1)}= \Ind _{\Om(w_2)}\g.
\end{equation}
This definition of $\phi(s_{w_2,w_1})$ is modeled on the rank one definition of $s_{u,v}$ given in Section \ref{Spielberg}. We refer the reader to \cite[Section 7]{rs2} for a proof of the isomorphism.

\begin{figure}[htbp]
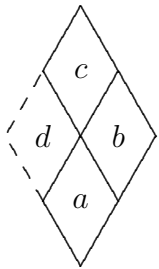
\label{H1}
\hfil
\centerline{
\beginpicture
\setcoordinatesystem units <0.5cm,0.866cm>  
\setplotarea x from -5 to 5, y from 0 to 4         
\put{$a$}   at  0 1
\put{$b$}   at  1  2
\put{$c$}   at  0 3
\put{$d$}   at  -1  2
\setlinear \plot -1 1  0 0   2 2  0 4  -1 3  0 2  1 3 /
\plot -1 1  0 2  1 1 /
\setdashes \plot  -1 1  -2 2  -1 3 /
\endpicture
}
\hfil
\caption{Tiles in an apartment.}
\end{figure}

A vital step in \cite[Section 7]{rs2} is the verification 
of the conditions (H0-H3) needed to define the higher rank Cuntz-Krieger algebra $\cA$ in Theorem \ref{mainboundaryaction}. Condition (H0) is obvious. Condition (H1) follows from the fact that in the configuration illustrated by Figure \ref{H1}, the tiles $a, b, c$ determine a unique tile $d$ lying in an apartment of $\cB$ containing $a, b, c$. Condition (H3) follows from thickness of the building, which allows words to be extended so as to lack periodicities.

\refstepcounter{picture}

The hardest condition to prove is (H2),  i.e. irreducibility of the associated directed graph.  This can be done by a direct combinatorial argument for $\tA_2$ groups. If the group $\G$ of Theorem \ref{mainboundaryaction} is a lattice subgroup of $\PGL_3(\KK)$, where $\KK$ is a local field of characteristic zero, this is done in \cite[Theorem 7.10]{rs2} using the Howe-Moore Ergodicity Theorem. In forthcoming work of T. Steger it is shown how to extend
the methods of the proof of the Howe--Moore Theorem  and so prove Theorem \ref{mainboundaryaction} in the stated generality.

\section{K-theory of rank $2$ Cuntz-Krieger algebras}\label{kr2}

According to the Classification Theorem for p.i.s.u.n. algebras, a higher rank Cuntz-Krieger algebra $\cA$ is classified by its K-groups and the class of the identity element in $K_0(\cA)$. It is therefore of some interest to compute the K-theory of these algebras. We have already observed that
for a fixed alphabet $A$ and fixed transition matrices $M_j$, the stable isomorphism class of  $\cA_D$ is independent of the decorating set $D$. For the purposes of computing $K_*(\cA)$, it is therefore enough to consider the algebra $\cA$ with trivial $D=A$. One can then follow the original approach to K-theory in the rank one case \cite{c3}, using the fact that $\cA \otimes\cK \cong \cF \rtimes \ZZ^r$, where $\cF$ is an AF algebra. A precise description of the K-theory was obtained in \cite{rs3} for the case $r=2$. For rank $r\ge 3$ the results are more complicated. We describe the rank $2$ results of \cite{rs3}.

Let $\cA$ be a rank two Cuntz-Krieger algebra associated with an alphabet $A$, trivial decorating set $D=A$ and matrices $M_1, M_2$, as in Section \ref{HRCK}. The matrices
$(I-M_1,I-M_2)$ and $(I-M_1^t,I-M_2^t)$ define homomorphisms $\ZZ^A\oplus \ZZ^A \to \ZZ^A$.
The K-theory of $\cA$ can be expressed as follows,
where $\tor(G)$ denotes the torsion part of a finitely generated abelian group $G$, and
$\rank(G)$ denotes the rank of $G$.

\begin{proposition}\label{Ktheory-}\cite[Proposition 2.14]{rs3} 
\begin{align*}
\rank(K_0(\cA))&=\rank(K_1(\cA))=
\rank(\coker\begin{smallmatrix}(I-M_1,&I-M_2)\end{smallmatrix})+
\rank(\coker\begin{smallmatrix}(I-M^t_1,&I-M^t_2)\end{smallmatrix}) \\
\tor(K_0(\cA))&\cong \tor(\coker\begin{smallmatrix}(I-M_1,&I-M_2)\end{smallmatrix}) \\
\tor(K_1(\cA))&\cong \tor(\coker\begin{smallmatrix}(I-M^t_1,&I-M^t_2)\end{smallmatrix}).
\end{align*}
\end{proposition}

If the algebras arise from group actions on buildings we can say more. 

\begin{theorem}\label{Ktheory}\cite[Theorem 3.2]{rs3} 
Let $\G$ be a group of automorphisms of a building $\cB$ of type $\tA_2$ which acts freely on the set of vertices of $\cB$ with finitely many orbits. Let $\Om$ be the boundary of $\cB$, so that $C(\Om)\rtimes \G$ is isomorphic to a higher rank Cuntz-Krieger algebra $\cA$.
If $M_1$,$M_2$ are the corresponding transition matrices, then
\begin{equation*}
K_0(\cA)=K_1(\cA)=\ZZ^{2n}\oplus \tor(\coker\begin{smallmatrix}(I-M_1,&I-M_2)\end{smallmatrix})
\end{equation*}
where $n=\rank(\coker\begin{smallmatrix}(I-M_1,&I-M_2)\end{smallmatrix})$.
\end{theorem}

\noindent The proof of Theorem \ref{Ktheory} uses symmetry considerations in the building.
Consider the special case where the action of $\G$ on the vertices of $\cB$ is also transitive, i.e. $\G$ is an $\tA_2$ group.  Extensive computational results are given in \cite{rs3} for more than 100 different groups with  $2\le q \le 11$, including all possible $\tA_2$ groups for $q=2,3$. In particular, for the group $\G= \G ({\rm C.1})$ of Example \ref{C1}, $K_0(C(\Om)\rtimes \G)=K_1(C(\Om)\rtimes \G)=(\ZZ/2\ZZ)^4\oplus (\ZZ/3\ZZ)$ and $[\id]=0$ in $K_0$.  This example is not typical in that $K_*$ usually has a free abelian component.

If $\G$ is an $\tA_2$ group then $[\id]$ is always a torsion element of $K_0(C(\Om)\rtimes \G)$. In fact \cite[Proposition 5.4]{rs3} proves that $(q^2-1)[\id]=0$.
Moreover, for $q \not\equiv 1 \pmod {3}$, $q-1$ divides the order of $[\id]$ and  for $q  \equiv 1 \pmod {3}$, $(q-1)/3$ divides the order of $[\id]$, \cite[Proposition 5.5]{rs3}. 
It follows that if $\G$ is an $\tA_2$ group and $q \ne 2,4$ then $[\id]$ is a nonzero torsion element in $K_0(C(\Om)\rtimes \G)$. Since $K_0=K_1$ and the $K_1$ group of a rank one Cuntz-Krieger algebra is torsion free, it is immediate that for $q \ne 2,4$, $C(\Om)\rtimes \G$ is not isomorphic to any rank one Cuntz-Krieger algebra. 

Abundant experimental evidence suggests that for algebras associated with $\tA_2$ groups it is always true that
$[\id]$ has order $q-1$ for $q \not\equiv 1 \pmod {3}$ and has order $(q-1)/3$ for $q  \equiv 1 \pmod {3}$.

In view of Theorem \ref{Ktheory}, it is worth considering the general structure of rank two Cuntz-Krieger algebras $\cA$ for which $K_0=K_1$ and for which this group has even rank. The following more general result applies.

\begin{proposition}
Let $\cA$ be a p.i.s.u.n. $C^*$-algebra with $K_0(\cA)=K_1(\cA)=\ZZ^{2n}\oplus T$ where $T$ is a finite abelian group.
Then $\cA$ is stably isomorphic to $\cA_1 \otimes \cA_2$,
where $\cA_1$, $\cA_2$ are simple rank one Cuntz-Krieger algebras.
\end{proposition}

\begin{proof}
By \cite[Proposition 6.6]{ror} we can find simple rank one Cuntz-Krieger algebras $\cA_1$, $\cA_2$ such that $K_*(\cA_1)=(\ZZ^{n}\oplus T,\ZZ^{n})$ and $K_*(\cA_2)=(\ZZ,\ZZ)$.
The K\"unneth Theorem for tensor products \cite[Theorem 23.1.3]{bl} shows that
$K_*(\cA_1 \otimes \cA_2)= (\ZZ^{2n}\oplus T,\ZZ^{2n}\oplus T$).
Since the algebras involved are all p.i.s.u.n. and satisfy the U.C.T., the result follows from the Classification Theorem \cite{k'}. 
\end{proof}

In particular, the algebras of Theorem \ref{Ktheory} are stably isomorphic to tensor products of rank one Cuntz-Krieger algebras.


\begin{thebibliography}{RRR}

\bibitem [An]{an} C. Anantharaman-Delaroche, Classification des $C\sp *$-alg\'ebres purement infinies nucl\'eaires (d'apr\`es E. Kirchberg). S\'eminaire Bourbaki, Vol. 1995/96. {\it Ast\'erisque} {\bf 241} (1997), 7--27.

\bibitem [Bl]{bl} B. Blackadar, {\it K-theory for Operator Algebras,
Second Edition}, MSRI Publications~5, Cambridge University Press,
Cambridge, 1998.

\bibitem [Br1]{bro1} K. Brown, {\it Buildings}, Springer-Verlag, New York, 1989.

\bibitem [Br2]{bro2} K. Brown, Five lectures on buildings, {\it Group Theory from a Geometrical
Viewpoint (Trieste 1990)}, 254--295, World Sci.\ Publishing, River Edge, N.J., 1991. 

\bibitem[Ca]{ca} D. I. Cartwright, A brief introduction to buildings,
{\it Harmonic Functions on Trees and Buildings (New York 1995)}, 45--77,
{\it Contemp.\ Math.} {\bf 206}, Amer.\ Math.\ Soc., 1997.

\bibitem[CMS]{cms}  D. I. Cartwright, W. M{\l}otkowski and T. Steger, Property (T) and
$\wt A_2$ groups, {\it Ann.\ Inst.\ Fourier} {\bf 44} (1993), 213--248.

\bibitem[CMSZ]{cmsz} D. I. Cartwright, A. M. Mantero, T. Steger and
A. Zappa, Groups acting simply transitively on the vertices of a
building of type~$\wt A_2$, I and II,\ {\it Geom.\ Ded.}  {\bf 47}
(1993), 143--166 and 167--223.

\bibitem[C1]{c1} J. Cuntz, Simple $C\sp*$-algebras generated by isometries, {\it Comm. Math. Phys.} {\bf 57} (1977), 173--185.

\bibitem[C2]{c2} J. Cuntz, A class of $C^*$-algebras and topological Markov chains:
Reducible chains and the Ext-functor for $C^*$-algebras, {\it Invent.\ Math.}
 {\bf 63} (1981), 23-50. 

\bibitem[C3]{c3} J. Cuntz, K-theory for certain $C^*$-algebras, {\it Ann.\ of Math.}
{\bf 113} (1981), 181-197.

\bibitem[CK]{ck} J. Cuntz and W. Krieger, A class of $C^*$-algebras and topological
Markov chains, {\it Invent.\ Math.} {\bf 56} (1980), 251-268.

\bibitem[DD]{dd} W. Dicks and M. J. Dunwoody, {\it Groups Acting on Graphs},  Cambridge University Press, Cambridge, 1989.  

\bibitem[G]{g} P. Garrett, {\it Buildings and Classical Groups}, Chapman \& Hall, London, 1997.

\bibitem [K1]{k} E. Kirchberg, Exact $C^*$-algebras, tensor products, and the classification
of purely infinite algebras, {\it Proceedings of the International
Congress of Mathematicians (Z\"urich, 1994)}, Vol.~2, 943--954,
Birkh\"auser, Basel, 1995.

\bibitem [K2]{k'} E. Kirchberg, The classification of purely infinite $C^*$-algebras
using Kasparov's theory, in {\it Lectures in Operator Algebras}, Fields Institute Monographs,
Amer. Math. Soc., 1998.

\bibitem[RR]{rr} J.~Ramagge and G.~Robertson, Triangle buildings and
actions of type $\text{\rm{III}}_{1/q^2}$, {\it J.\ Func.\ Anal.} {\bf 140} (1996),
472--504.

\bibitem[RRS]{rrs} J.~Ramagge, G.~Robertson and T. Steger, A Haagerup inequality for $\tilde A_1 \times\tilde A_1$ and $\tilde A_2$ buildings,
{\it Geometric and Funct. Anal.} {\bf 8} (1998), 702--731.

\bibitem[RS]{rs} G. Robertson and T. Steger, Maximal abelian subalgebras of the
group factor of an $\widetilde A_2$ group, {\it J. Operator Theory } {\bf
36} (1996), 317--334.

\bibitem[RS1]{rs1} G. Robertson and T. Steger,  $C^*$-algebras arising from group actions on the boundary of a triangle building, {\it Proc.\ London Math.\ Soc.}  {\bf 72} (1996), 613--637.

\bibitem[RS2]{rs2} G. Robertson and T. Steger, Affine buildings, tiling systems and higher rank Cuntz-Krieger algebras, {\it J. reine angew. Math.} {\bf 513} (1999), 115--144.

\bibitem[RS3]{rs3} G. Robertson and T. Steger, K-theory for rank two Cuntz-Krieger algebras, {\it preprint}.

\bibitem[Ron]{ron} M. Ronan, {\it Lectures on Buildings}, Perspectives
in Mathematics, Vol.~7,
Academic Press, New York, 1989.

\bibitem[Ror]{ror} M. R\o rdam,  Classification of Cuntz-Krieger algebras, {\it K-theory}  {\bf 9} (1995), 31--58.

\bibitem[Ser]{ser} J.-P. Serre, {\it Trees}, Springer-Verlag, Berlin, 1980.

\bibitem[Sp]{s} J. Spielberg, Free product groups, Cuntz-Krieger algebras, and covariant maps,
{\it International J. Math.} {\bf 2} (1991), 457-476.

\bibitem[St]{st} T. Steger, Local fields and buildings,
{\it Harmonic Functions on Trees and Buildings (New York 1995)}, 79--107,
Contemp.\ Math. {\bf 206}, Amer.\ Math.\ Soc., 1997.

\end{thebibliography}
\end{document}